\documentstyle[12pt]{article}
\title{ On Disjointness  of Mixing Rank One Actions}
\author{V.V. Ryzhikov}
\date{}
\def\N{{\bf N}}
\def\eps{\varepsilon}

\def\u{\bigsqcup}
\def\eps{\varepsilon}
\date{vryzh@mail.ru}
\begin{document}

\maketitle
\begin{abstract}  
For flows the rank is an invariant by linear change of time.  But what we can say about isomorphisms? It seems that in case of mixing flows this problem is the most difficult. However the known technique of joinings provides non-isomorphism for mixing rank-one flows  under linear change of time. For automorphisms we consider another problems (with similar solutions). For example, the staircase cutting-and-stacking construction is determined  by  height $h_1$ of the first tower and a sequence $\{r_j\}$  of cut numbers. Let us consider two  similar constructions: one is set by $(h_1, \{r_j\})$, another is set by ($h_1+1, \{r_j\}$), and $r_j=j$. We prove a general  theorem implying  the non-isomorphism of these constructions.
\end{abstract}
\vspace{3mm}
{\large

\section{Introduction.  Rank-one transformations and flows. Joinings and Disjointness. }
S. Kalikow has proved \cite{K} that mixing rank-one transformations  are 3-fold mixing.  In \cite{R}   we have got a  joining proof of Kalikow's theorem (see also \cite{R2}) via ``powder'' method: a pairwise independent joining $\nu$ has to have a powder part ($Di(\nu)>0$)  that trivializes the joining ($\nu$ must be a product measure).   Now we ``project'' this method  into a two-dimensional  situation to trivialize  joinings  of two rank-one transformations. We present Theorem 1 on the disjointness for 
certain pairs of rank one mixing transformations of a probability  space $(X,\mu)$.  For example,  let $T$ be a staircase transformation defined by
a parameter  $h_1$ (the  height  of "the first tower")  and a cutting-sequence $r_j=j$, see \cite{Ad}.
Let $\tilde T$ be a staircase transformation with   $\tilde h_1 >h_1$ and $\tilde r_j=j$.  Then $T$ and $\tilde T$  are disjoint.  

We recall that   
 El H. El Abdalaoui  developed Bourgain's method and  showed that Ornstein's stochastic transformations $T$, $\tilde T$ are almost surely (spectrally) disjoint \cite{A}. We present  sufficient conditions for the disjointness of two rank-one transformations. This gives  a large class of  pairwise non-isomorphic explicit examples.

In \cite{R1} we  stated that a mixing rank-one flow $T_t$ was disjoint with any flow $T_{\alpha t}$ as $\alpha >1$. We present a poof now (Theorem 2).  In connection with  preprint  \cite{P} let's note that  for a (mixing rank-one) flow $T_t$ with Lebesgue spectrum for all  $\alpha \neq 0$ the flows  $T_{\alpha t}$  are  spectrally isomorphic.

\bf Definitions. \rm An automorphism\footnote{a measure-preserving invertible transformation of a probability Lebesgue space $(X,\mu)$}
$T:X\to X$ is said to be of
{\it rank one}, if  there is a sequence $\xi_j$ of measurable partitions
of $X$ in the form
$$\xi_j= \{ E_j,\    TE_j, \  T^2 E_j,\ \dots,   T^{h_j}E_j,  E^\prime_j\}$$
 such that the partitions $\xi_j$ converge to the partition onto points ($\xi_j \to\eps$).
The collection $$E_j, TE_j, \  T^2 E_j,\ \dots,   T^{h_j}E_j$$ is called
Rokhlin's tower ( ${E}_j^\prime=X\setminus \u_{i=0}^{h_j} T^i E_j$).

 The property ``to be   rank-one flow'' is defined as the existence of a continuous rectangle-tower sequence $\xi_j$ such that any measurable
set $A$ can be approximated by a $\xi_j$-measurable set $A_j$ ( a measurable union  of floors in $j$-tower).  We will use only the following property of rank-one flows:  
\it for some sequences $t_j\to 0$, $h_j\to\infty$  there is a  sequence $E_j$  such that   \rm
$$\xi_j= \{ E_j,\    T_{t_j}E_j, \  T_{2 t_j} E_j,\ \dots,   T_{h_j t_j}E_j,  E^\prime _j\}\to \eps.$$
We suppose without  a loss of generality that $1/t_j\in \N$. For a flow $\tilde T_t=T_{\alpha t}$ we  find easy $\tilde E_j $  setting  $\tilde t_j=\alpha t_j$ and 
 $\tilde h_j=[h_j/\alpha]$.  So $\tilde T_t$ is a  rank-one flow as well. 

A {\it joining} of automorphisms $\tilde T$ and $T$  is defined to be a $(\tilde T\times T)$-invariant
measure $\nu$ on $X\times X$ with its marginals  equal to  $\mu$:
$$\nu(A\times X)=\nu(X\times A)=\mu(A).$$
A joining $\nu$ is called ergodic if the dynamical system
$(\tilde T\times T, X\times X, \nu)$ is ergodic.

If $\mu\times\mu$ is a unique joining of $\tilde T$ and $T$, then (obviously non-isomorphic) $\tilde T$ and $T$ are called \it disjoint,  \rm   see \cite{F}.

We say that $T$ is {\it mixing} (or 2-fold mixing) if
  for all measurable sets $A,B$
$$
\mu(A\cap T^mB) \to  \mu(A)\mu(B),   \ m\to\infty. $$
\vspace{3mm}

\section{Auxiliary assertions}

\bf LEMMA 1.1. (Blum-Hanson). \it Let  a
 sequence $\lbrace a^{z}_j \rbrace , z,j\,\in\,\bf N$,
satisfy the conditions:
$$\sum_z a^{z}_j =1,\  \ a^{z}_j  \geq 0, \ \ and \ 
  max_z \lbrace a^{z}_j \rbrace
\,\to \,0,\ \  as \ \ j\,\to \,\infty.$$
If $T$ is mixing, then 
$$
  \left\| \sum_z a^{z}_j T^zf - \int \! f\, \right\|_2 \,\to \,0.
$$
 \rm
\vspace{3mm}

Proof.  Let $\int  f=0$.  Put $P_j = \sum_z a^{z}_j T^z.$
Let us show that $  \| P_jf\|_2 \,\to \,0.$
One has
$$P_j^\ast P_j = \sum_w b^{w}_j T^w ,$$
where  the sequence $\lbrace b^{w}_j \rbrace$ satisfies
$$   b_{w}^j \leq \sum_z a^{w-z}_j a^{z}_j \leq max_z {a^z_j } \to \,0.
$$
Since $T$ is mixing, one has   $\sum_w b^{w}_j T^w f \to 0$ (weakly).
Thus,
$$\| P_j f\|^2=(P_j^\ast P_j f, f)\to 0.$$

If $\int  f\neq 0$, we get
$$ \|P_jf-\int f\|^2\to 0.$$
$\Box$
\vspace{3mm}

\bf LEMMA 1.2. \it Let 
 a sequence $\lbrace a^{z}_j \rbrace , z\,\in\N$,
satisfy the conditions:
$$
\sum_z a^{z}_j =1, \ a^{z}_j  \geq 0; \ \ max_{z} \lbrace \sum_{w=z}^{z+1/t_j} a^{w}_j \rbrace
\,\to \,0,\,j\,\to \,\infty.$$
If a flow $T_t $ is mixing, then 
$$
  \left\| \sum_z a^{z}_j T_{zt_j}f - \int \! f\, \right\|_2 \,\to \,0.
$$
 \rm
\vspace{3mm}
\\
\bf  LEMMA  2.1.  \it Let $T$ be a rank-one transformation with a corresponding sequence $E_j$. We set 
 $a^{z}_j =\mu(T^zE_j | E_j)$.
If $T$ is mixing,  then  $\lim_j max_{z>0} \{a^z_j\} = 0$.  \rm
\vspace{3mm}

Proof. We have  $ max_{z>0} \{a^z_j\}= max_{z>h_j} \{a^z_j\}$.
Suppose $\lim_j max \{a^z_j :z>h_j\} = a>0$,
 $\mu(T^{z_j}E_j | E_j)\to a,$ hence,
$$\mu(T^{z_j}T^kE_j | T^kE_j)\to a, \ \ (0\leq k\leq h_j).$$
We can approximate the measurable set $A$  by $\xi_j$-measurable
sets $A_j$ ($A_j$ is a union of certain floors $T^kE_j$). We have
$\limsup_j \mu(T^{z_j}A_j | A_j)\geq  a,$  hence,  for all $A$, $\mu(A)>0$,  
$\limsup_j \mu(T^{z_j}A | A)\geq  a$
holds. The mixing implies 
$\mu(T^{z_j}A | A)\to \mu(A).$  
Thus,  $a\leq \mu(A), \ \ a=0.$
$\Box$
\\
\bf  LEMMA  2.2  (On  little by little  Returning). \it Let $T_t$ be a mixing rank-one flow. Then setting   $a^{w}_j =\mu(E_j | T_{w t_j} E_j )$
 we have $\lim_j max_{ z >0} \sum_{w=z}^{z+\frac{1}{t_j}} a^w_j = 0.$  \rm
\vspace{3mm}
\\
Let's denote $E_j^1=\bigcup_{w= 0}^{1/t_j}   T_{wt_j}E_j$.  From  Lemma 2.2 we see  that  
 $\mu(E_j^1 | T_{z_j t_j} E_j )\to 0$  for any  sequence  $z_j\to+\infty$. 

\section{Disjointness of Transformations}            
\bf THEOREM 1.   \it Let $\tilde T,T$ be  rank-one transformations with    height sequences 
$\tilde  h_j$ and $h_j$, respectively.
Let $\nu$ be an ergodic joining of $\tilde T$ and $T.$ If $T$ (or $\tilde T$) is mixing, and $${\tilde h_j}/h_j\to \alpha \in (0,1),$$
then $\nu=\mu\otimes\mu$, i.e. $\tilde T$ and $T$ are disjoint.  \rm
\vspace{3mm}

{Theorem 1 has been  presented   
at {\it Laboratoire de Mathematiques Raphael Salem }  of  Rouen University. 
The author  thanks  El H. El Abdalaoui, T.~de~la~Rue  and J.-P. Thouvenot for  discussions.}

Proof.
For  $\eps >0$ let us  define a set $ D_{\eps, j}$ of $\eps$-light block  indexes:
$$
  D_{\eps, j}=\{z\in [0,\tilde h_j]\times [0, h_j] \ :\  \nu(V^z_j)<\eps \mu(E_j)\}, $$
 where $z=(z1,z2)$,
$V^z_j= \tilde T^{z1}\tilde E_j\times  T^{z2} E_j$.
Now we calculate the total mass $Di(\nu)$ of infinitely light blocks, i.e.  an asymptotically  diffused portion (a powder) of a joining.
$$
  Di(\nu)=\lim_{\eps\to 0}
\left(\limsup_{j\to \infty} \sum_{z\in D_{\eps,j}}\nu(\bar{E}^z_j)\right).
$$
   Now we show how heavy blocks generate  light blocks.

\bf LEMMA 3.  \it
If $\nu$ is a joining of $\tilde T$ and $T$, and ${\tilde h_j}/h_j\to \alpha \in (0,1),$
then $Di(\nu)>0$ (there is a powder).  \rm
\vspace{3mm}

Proof.   Lemma 2.1 and the following picture  show   that  heavy blocks
under the action of  some powers of $\tilde T\times T$  generate
many light blocks, hence,  $Di(\nu)>0$ forever.

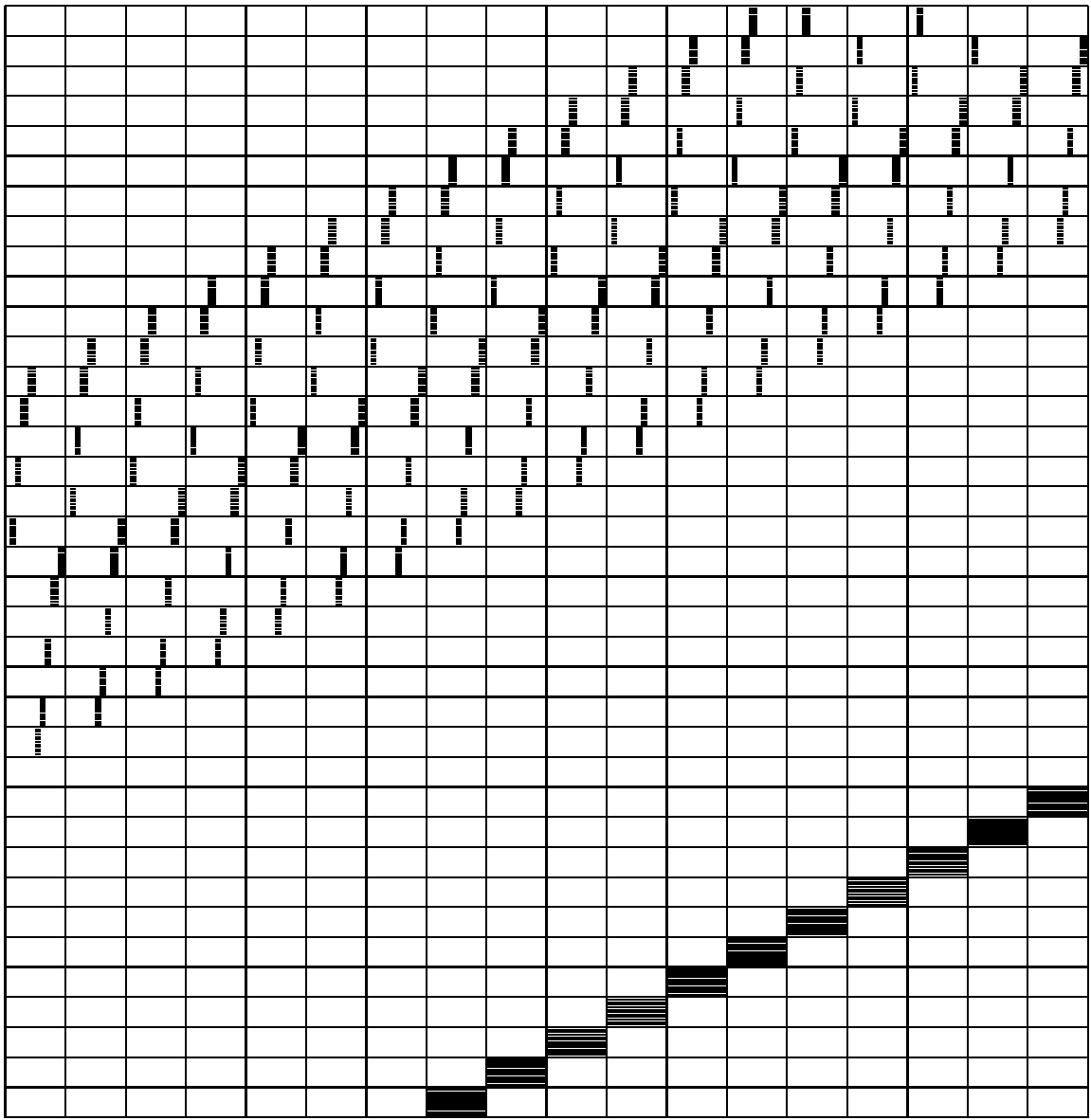
\begin{picture}(0, 450) 
\multiput(0,0)(0,12){38}%
{  \line(1, 0){432}   }

\multiput(0,0)(24,0){19}%
{  \line(0, 1){444}   }
\multiput(84,0)(24,12){11}%
{\multiput(84,0)(0,1){12}
{  \line(1, 0){24}   }}

\multiput(12,144)(24,12){18}%
{\multiput(0,0)(0,1){12}
{  \line(1, 0){2}   }}

\multiput(14,156)(24,12){18}%
{\multiput(0,0)(0,1){12}
{  \line(1, 0){2}   }}
\multiput(16,180)(24,12){18}%
{\multiput(0,0)(0,1){12}
{  \line(1, 0){2}   }}
\multiput(18,204)(24,12){18}%
{\multiput(0,0)(0,1){12}
{  \line(1, 0){3}   }}
\multiput(21,216)(24,12){18}%
{\multiput(0,0)(0,1){12}
{  \line(1, 0){3}   }}

\multiput(0,0)(-12,72){2}
{
\multiput(14,156)(24,12){17}%
{\multiput(0,0)(0,1){12}
{  \line(1, 0){2}   }}
\multiput(16,180)(24,12){16}%
{\multiput(0,0)(0,1){12}
{  \line(1, 0){2}   }}
\multiput(18,204)(24,12){14}%
{\multiput(0,0)(0,1){12}
{  \line(1, 0){3}   }}
\multiput(21,216)(24,12){13}%
{\multiput(0,0)(0,1){12}
{  \line(1, 0){3}   }}
}

\end{picture}  

\begin{center}
$Di(\nu)=0$ on the right below  implies $Di(\nu)>0$ at the top.
\end{center}
\vspace{2mm}

{\bf LEMMA 4.} \it
If 
 $Di(\nu)>0$, then $\nu =\mu\otimes \mu$.
\rm

{\large
Proof. 
Let's show how light blocks trivialize a joining.   The Blum-Hanson lemma and rank one approximations will be our tools.              
We define  columns  in the following way:
$$
C_j^{w}=\u_{i=0}^{\delta h_j}
\tilde T^{w+i}E_j\times T^{i}E_j.
$$
Given small $\delta >0$ we  find  a sequence of
sets  $F_j$ of the form
$$
F_j =:\u_{h\in D_j} (Id\times  T^h)C_j, \ \ C_j= C_j^{w_j},$$
for some sequences
 $$D_j\subset \{0,1, \dots (1-\delta) h_j\},
\ \ w_j\in \{0,1, \dots (1-\delta) h_j\}.$$

  From   $Di(\nu)>0$ it follows that there is a sequence of 
 $F_j$
with   $D_j$  that numerate  $\eps_j$-light columns only  ($\nu((Id\times  T^h)C_j)\leq \eps_j\to 0$), and 
  $$\nu(F_j)\to a>0.$$
The sets $D_j$ will  satisfy the condition:
$$\max_{h\in D_j}\{a_j^h\} \to 0, \  j\to\infty, \ \ \
 \sum_{h\in D_j}a_j^h =1,
$$
where
$a_j^h=\nu((Id\times  T^h)C_j   \,|\, F_j).$
  Since $F_j$  are almost invariant
 with respect to $\tilde T\times T$   and
 $\nu$ is ergodic with respect to $\tilde T\times T$,  we  get 
$$\nu(\ | F_j)\to \nu.$$
Let's show 
$$\nu(\ | F_j)\to \ \mu\times\mu.$$

If $A,B$ are $\xi_k$-measurable, then for all $j\geq k$ 
$$\nu(A\times B\ | F_j)=\sum_{h\in D_j}a_j^h \nu(A\times T^{-h}B\ | C_j)=
\sum_{h\in D_j}a_j^h \lambda(A\times T^{-h}B\ | C_j),$$  where  $\lambda=\mu\times\mu$. 
 Lemma 1.1 provides 
$$\sum_{h\in D_j} a_j^h T^{-h}\chi_{ B} \ \to_{L_2}  \ Const \
\equiv \ \mu(B) \   (j\to\infty).
$$
Thus,
$$\nu(A\times B) =
\lim_{j\to\infty} \nu(A\times B  |\, F_j) =\lim_{j\to\infty}\sum_{h\in D_j}a_j^h \nu(A\times T^{-h}B\ | C_j)
$$
$$
=\lim_{j\to\infty} \ 
\int_{X\times X} \chi_A   \otimes 
\left( \sum_{h\in D_j} a_j^h T^{-h}   \chi_{B}\right)\, d\lambda( \ |C_j)=
\mu(B)\,\nu(A\times  X) =
\mu(A)\mu(B).
$$
We used above the fact that marginal projections of $\lambda( \ |C_j)$ were of  densities less than $\delta^{-1}$.   Finally,
$\nu= \lim \nu(\ | F_j)= \mu\times\mu.$
$\Box$     

Theorem 1 follows from   Lemma 4.  (We use now the mixing for both $\tilde T$ and $T$. In fact, it's enough for a proof that one of them is  mixing. Let it be an exercise.)  
}

\section{Disjointness of Flows}
\bf THEOREM 2. \cite{R1} \it  Let a  rank-one flow $T_t$ be mixing.  Then it is disjoint from any flow $T_{\alpha t}$ 
as $\alpha >1$. \rm
\\
Proof.  In fact we repeat the proof of Theorem 1 with little modifications that we try to explain below.  Let's look, for example,  at a joining $\nu$ of  $T_{2t}$ and $T_t$.

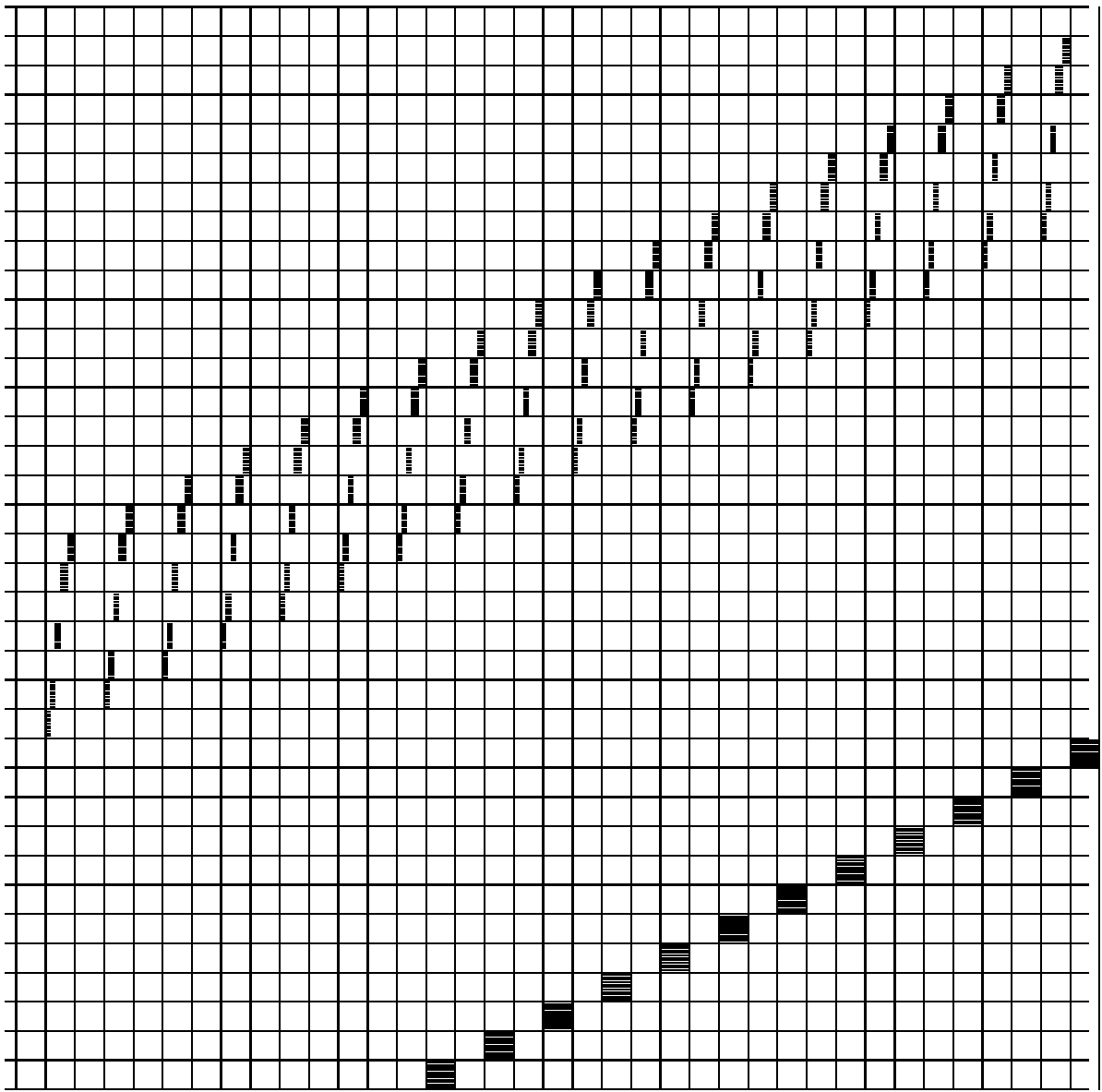
\begin{picture}(0, 450) 
\multiput(0,0)(0,12){38}%
 {\line(1, 0){444}}   

\multiput(0,0)(12,0){38}%
{  \line(0, 1){444}   }


\
\multiput(84,0)(24,12){12}%
{\multiput(84,0)(0,1){12}
{  \line(1, 0){12}   }}

\multiput(12,144)(24,12){18}%
{\multiput(0,0)(0,1){12}
{  \line(1, 0){2}   }}

\multiput(14,156)(24,12){18}%
{\multiput(0,0)(0,1){12}
{  \line(1, 0){2}   }}
\multiput(16,180)(24,12){18}%
{\multiput(0,0)(0,1){12}
{  \line(1, 0){2}   }}
\multiput(18,204)(24,12){18}%
{\multiput(0,0)(0,1){12}
{  \line(1, 0){3}   }}
\multiput(21,216)(24,12){18}%
{\multiput(0,0)(0,1){12}
{  \line(1, 0){3}   }}
\end{picture} 
\begin{center}
Images of $T_{w}E_j^1\times E_j^1$ under the action  of $T_{2z}\times T_{z}$.
\end{center}

If near  the right bottom corner a joining $\nu$ is concentrated in a union 
$$Band_{(w_j,0)} = \bigcup_{h=0}^ { {1}/{t_j}} \bigcup_{z=0}^ { h_j -w_j} T_{(2z+w_j)t_j}E_j\times T_{(z+h)t_j}E_j,$$
then left above  we observe  its dispersion (applying Lemma 1.2).    
This dispersion means that an essential  part of $\nu$ is situated in a big collection of ``left bands''
as $$Band_{(0,v)}
= \bigcup_{h=0}^ { {1}/{t_j}} \bigcup_{z=0}^ { w_j} T_{2zt_j}E_j\times T_{(z+h+v)t_j}E_j,$$  and each of these bands has  small $\nu$-measure. Thus, we have $Di(\nu)>0$ with respect to partitions 
$$\xi_j^1= \{ E_j^1,\    T_{1}E_j^1, \  T_{2} E_j^1,\ \dots,   T_{H_j}E_j^1\},$$
where $h_j t_j/H_j\approx 1$,   $E^1_j=\u_{w=  0}^{1/t_j}   T_{wt_j}E_j$.
 (Warning: $\xi^1_j\to\eps$ fails. So, we have to deal with $\xi_j$-approximations.)
 Arguing  as in the proof of a lemma 4 by use of Lemma 2.2     we get  $\nu=a\mu\times\mu+\dots$ for some $a>0$.   Assuming $\nu$ to be ergodic  we get  $\nu=\mu\times\mu$. Then we remember  that every joining
is a convex   sum (integral) of ergodic ones.  $\Box$}
\vspace{5mm}

Remark.  There is a possibility to construct  two  rank-one mixing explicit transformations with disjoint spectra.  
J.-P. Thouvenot  asks: is  the rank invariant with respect to the  spectral isomorphism?  
Disjointness and Spectral  Disjointness,  whether these two concepts coincide for rank-one
(mixing) transformations?

 \end{document}